\documentclass[11pt]{amsart}
\usepackage{amsmath, amsthm, latexsym, amssymb, amsfonts, epsfig, epsf, color}

\newenvironment{pf}{\proof[\proofname]}{\endproof}
\theoremstyle{plain}
\newtheorem{Th}{Theorem}[section]
\newtheorem{Cor}[Th]{Corollary}

\newtheorem{Prop}[Th]{Proposition}
\newtheorem{Lemma}[Th]{Lemma}
\numberwithin{equation}{section} \theoremstyle{definition}
\newtheorem{Rem}[Th]{Remark}

\newtheorem{Ex}[Th]{Example}

\newtheorem{Def}[Th]{Definition}

\newcommand{\cal}[1]{\mathcal{#1}}

\newcommand{\C}{\mathbb C}
\newcommand{\K}{\mathbb K}
\newcommand{\mS}{\mathbb S}
\newcommand{\Z}{\mathbb Z}
\newcommand{\R}{\mathbb R}

\newcommand{\D}{\Delta}
\newcommand{\cA}{\cal A}

\newcommand{\cF}{\cal F}
\newcommand{\cC}{\cal C}

\newcommand{\conv}{\operatorname{conv}}

\newcommand{\Vol}{\operatorname{Vol}}
\newcommand{\la}{\langle}
\newcommand{\ra}{\rangle}

\newcommand{\rs}[1]{Section~\ref{S:#1}}
\newcommand{\rl}[1]{Lemma~\ref{L:#1}}
\newcommand{\rp}[1]{Proposition~\ref{P:#1}}

\newcommand{\rr}[1]{Remark~\ref{R:#1}}
\newcommand{\rex}[1]{Example~\ref{Ex:#1}}
\newcommand{\re}[1]{(\ref{E:#1})}

\newcommand{\rc}[1]{Corollary~\ref{C:#1}}
\newcommand{\rt}[1] {Theorem~\ref{T:#1}}

\newcommand{\rf}[1]{Figure~\ref{F:#1}}

\begin{document}

\title[Criteria for strict monotonicity of the mixed volume]{Criteria for strict monotonicity of the mixed volume of convex polytopes}
\author{Fr\'ed\'eric Bihan}
\address[Fr\'ed\'eric Bihan]{D\'epartement de Math\'ematiques LAMA, 
Universit\'e Savoie Mont Blanc, France}
\email{frederic.bihan@univ-smb.fr}
\author{Ivan Soprunov}
\address[Ivan Soprunov]{Department of Mathematics\\ Cleveland State University\\ Cleveland, OH USA}
\email{i.soprunov@csuohio.edu}

\begin{abstract} Let $P_1,\dots, P_n$ and $Q_1,\dots, Q_n$ be convex polytopes in $\R^n$ such that $P_i\subset Q_i$.
It is well-known that the mixed volume has the monotonicity property: $V(P_1,\dots,P_n)\leq V(Q_1,\dots,Q_n)$.
We give two criteria for when this inequality is strict
in terms of essential collections of faces as well as mixed polyhedral subdivisions.
This geometric result allows us to  characterize sparse polynomial systems with Newton polytopes $P_1,\dots,P_n$ whose number of isolated solutions equals the normalized volume of the convex hull of $P_1\cup\dots\cup P_n$.
In addition, we obtain an analog of Cramer's rule for sparse polynomial systems.
\end{abstract}

\maketitle

\section{Introduction}
The mixed volume is one of the fundamental notions in the theory of convex bodies. 
It plays a
central role in the Brunn--Minkowski theory and in the theory of sparse polynomial systems.
The mixed volume is the polarization of the volume form
on the space of convex bodies in $\R^n$.
More precisely, let $K_1,\dots, K_n$ be $n$ convex bodies in $\R^n$ and 
$\Vol_n(K)$ the Euclidean volume of a body $K\subset \R^n$. Then the mixed volume of $K_1,\dots,K_n$
is
\begin{equation}\label{E:MV}
V(K_1,\dots,K_n)=\frac{1}{n!}\sum_{m=1}^n(-1)^{n+m}\!\sum_{i_1<\dots<i_m}\Vol_n(K_{i_1}+\dots+K_{i_m}),
\end{equation}
where $K+L=\{x+y\in\R^n \ |\ x\in K, y\in L\}$ denotes the Minkowski sum of bodies $K$ and $L$.
It is not hard to see from the definition that the mixed volume is 
symmetric, multilinear with respect to Minkowski addition, invariant under translations, and coincides with the volume on the diagonal,
i.e. $V(K,\dots,K)=\Vol_n(K)$. What is not apparent from the definition is that
it satisfies the following {\it monotonicity property}, see \cite[(5.25)]{S}.
 If $L_1,\dots, L_n$ are convex bodies such that $K_i\subseteq L_i$
for $1\leq i\leq n$ then
$$V(K_1,\dots,K_n)\leq V(L_1,\dots,L_n).$$
The main goal of this paper is to give a geometric criterion for strict monotonicity in the class of convex polytopes.
We give two equivalent criteria in terms of essential collections of faces and Cayley polytopes, see \rt{main2} and \rt{main3}. The first criterion is especially simple when all $L_i$ are equal (\rc{MV=V}) which is the situation in our application to sparse polynomial systems.
In the general case of convex bodies this is still an open problem, see \cite[pp. 429--431]{S} for special cases and conjectures. 

The role of mixed volumes in algebraic geometry originates in the work of Bernstein, Kushnirenko, and Khovanskii, who gave a vast generalization of the classical B\'ezout formula for the intersection number of hypersurfaces in the projective space, see \cite{Be, Kho, Kush}. This beautiful result, which links algebraic geometry and convex geometry through toric varieties and sparse polynomial systems, is commonly known as the BKK bound. Consider an $n$-variate Laurent polynomial system $f_1(x)=\cdots=f_n(x)=0$ over an algebraically closed field $\K$. The {\it support} $\cA_i$ of $f_i$ is
the set of exponent vectors in $\Z^n$ of the monomials appearing with a non-zero coefficient in $f_i$. The {\it Newton polytope} $P_i$ of $f_i$
is the convex hull of $\cA_i$.
The BKK bound states that the number of isolated solutions of the system in the algebraic torus $(\K^*)^n=(\K\setminus\{0\})^n$  is at most $n! \, V(P_1,\dots,P_n)$.
Systems that attain this bound must satisfy a {\it non-degeneracy condition}, which means that certain
subsystems (predicting solutions ``at infinity") have to be inconsistent, see \rt{BKK}. 
However, the non-degeneracy condition may be hard to check.
Let $\cA=\cup_{i=1}^n\cA_i$ be the total support of the system and choose an order of its elements, $\cA=\{a_1,\ldots,a_{\ell}\}$. Then the system
can be written in a matrix form
\begin{equation}\label{matricial}
C \,x^{A}=0,
\end{equation}
where $C\in \K^{n \times \ell}$ is the matrix of coefficients, $A \in \Z^{n \times \ell}$ is the matrix of exponents whose
columns are $a_1,\ldots,a_{\ell}$, and $x^{A}$ is the transpose of
$(x^{a_1}, \dots, x^{a_{\ell}})$, see \rs{pol}. The solution set of \eqref{matricial} in $(\K^*)^n$ does not change after left multiplication of $C$ by a matrix in $\mbox{GL}_{n}(\K)$. Such an operation does not preserve the individual supports of \eqref{matricial} in general, but preserves the total support $\cA$, see \rr{left}. Furthermore, let $\bar{A} \in \Z^{(n+1) \times \ell}$ be the {\it augmented exponent matrix}, obtained by appending a first row of $1$ to $A$. Then left multiplication of $\bar{A}$ by a matrix in $\mbox{GL}_{n+1}(\Z)$ with first row $(1,0,\dots, 0)$ 
corresponds to a monomial change of coordinates of the torus $(\K^*)^n$ and a translation of $\cA$, hence, does not change the number of solutions of the system in $(\K^*)^n$, see \rs{pol}.


Assume all  Newton polytopes of a system \eqref{matricial} are equal to some polytope $Q$. Then the number of isolated solutions of \eqref{matricial} is at most $n! V(Q,\ldots,Q)=n! \Vol_n(Q)$, by the BKK bound. In \rs{pol} we characterize systems that reach (or do not reach) this bound. This characterization follows from the geometric criterion of \rc{MV=V}, but it has a natural interpretation
in terms of the coefficient matrix $C$ and the augmented exponent matrix $\bar A$, see \rt{Ber}. In particular, it says that if $Q$ has a proper face such that the rank of the corresponding submatrix of $C$ (obtained by selecting the columns indexed by points of $\cA$ which belong to that face)
is strictly less than the rank of the  corresponding submatrix of $\bar{A}$, then \eqref{matricial} has strictly less than $n! \Vol_n(Q)$ isolated solutions in $(\K^*)^n$. Naturally, this characterization is invariant under the actions on $C$ and $\bar A$ described above.

Another consequence of \rt{Ber} can be thought of as a generalization of Cramer's rule for linear systems.
Linear systems occur when each $P_i$ is contained in the standard unit simplex $\Delta$ in $\R^n$.
The BKK bound for systems with all Newton polytopes equal to $\Delta$ is just $1=n!\, \Vol_n(\Delta)$ and, by Cramer's rule, if all maximal minors of the coefficient matrix $C$ are non-zero, then the system \eqref{matricial} has precisely one
solution in $(\K^*)^n$. We generalize this to an arbitrary Newton polytope $Q$: If no maximal minor of $C$ vanishes then the system \eqref{matricial} has the maximal number $n! \Vol_n(Q)$ of isolated solutions in $(\K^*)^n$, see Corollary \ref{C:nice}.

\subsection*{Acknowledgement} This project began at the Einstein Workshop on Lattice Polytopes at
Freie Universit\"at Berlin in December 2016. We are grateful to M\'onica Blanco, Christian Haase,
Benjamin Nill, and Francisco Santos for organizing this wonderful event and to the Harnack Haus
for their hospitality. We thank Gennadiy Averkov and Christian Haase for fruitful discussions. Finally, we are thankful to the anonymous referee for their comments and suggestions that led to a substantial improvement of the exposition.

\section{Preliminaries}\label{Pre}
In this section we recall necessary definitions and results from convex geometry and set up notation.
In addition, we recall the notion of essential collections of polytopes for which we give several equivalent definitions,
as well as define mixed polyhedral subdivisions and the combinatorial Cayley trick.

Throughout the paper we use $[n]$ to denote the set $\{1,\dots,n\}$.

\subsection*{Mixed Volume} For a convex body $K$ in $\R^n$ the function $h_K:\R^n\to\R$, 
given by $h_K(u)=\max\{\la u, x\ra\ |\ x \in K\}$ is the {\it support function} of $K$. Here  $\la u, x\ra$ is the standard scalar product in $\R^n$. 
For every $u \in\R^n$, we write $H_K(u)$ to denote the
{\it supporting hyperplane} for $K$ with outer normal $u$
$$H_K(u)=\{x\in\R^n \ |\ \la u,x\ra = h_K(u)\}.$$ 
Throughout the paper we use
$$K^u=K\cap H_K(u)$$ 
to denote the {\it face} of $K$ corresponding to the supporting hyperplane $H_K(u)$. 
Since $H_K(u)$ and $K^u$ are invariant under rescaling $u$ by a non-zero scalar, we often assume that 
when $u\neq 0$, it lies in the unit sphere $\mS^{n-1}\subset\R^n$. 
Clearly, for $u=0$ we have $H_K(u)=\R^n$ and $K^u=K$.

Let $V(K_1,\dots,K_n)$ be the $n$-dimensional mixed volume of
$n$ convex bodies $K_1,\dots, K_n$ in $\R^n$, see \re{MV}. 
We have the following equivalent characterization.
\begin{Th}\cite[Theorem 5.1.7]{S}\label{T:MV}
Let $\lambda_1,\dots, \lambda_n$ be non-negative real numbers. Then 
$\Vol_n(\lambda_1K_1+\dots+\lambda_nK_n)$ is a polynomial in $\lambda_1,\dots, \lambda_n$ whose
coefficient of the monomial $\lambda_1\cdots\lambda_n$ equals $V(K_1,\dots,K_n)$.
\end{Th}

\subsection*{Essential Collections} Throughout the paper we use ``collection'' as a synonym for ``multiset".
 Let $K_1,\dots, K_m$ be convex bodies in $\R^n$, not necessarily distinct.
We say that a collection $\{K_1,\dots, K_m\}$ is  {\it essential} if for
any subset $I\subset [m]$ of size at most $n$ we have

\begin{equation}\label{E:essential}
\dim\sum_{i\in I}K_i\geq |I|.
\end{equation}

Note that every sub-collection of an essential collection is essential. Also $\{K,\dots,K\}$, where $K$ is repeated
$m$ times, is essential if and only if $\dim K\geq m$.

The following well-known result asserts that essential collections of $n$ convex bodies characterize positivity of  the mixed volume.

\begin{Th}\cite[Theorem 5.1.8]{S}\label{T:essential} Let $K_1,\dots, K_n$ be $n$ convex bodies in $\R^n$. 
The following are equivalent:
\begin{enumerate}
\item $V(K_1,\dots,K_n)>0$;
\item  There exist segments $E_i\subset K_i$ for $1\leq i\leq n$ with linearly independent directions;
\item $\{K_1,\dots, K_n\}$ is essential.
\end{enumerate}
\end{Th}

Another useful result is the inductive formula for the mixed volume, see \cite[Theorem 5.1.7, (5.19)]{S}.
We present a variation of this formula for polytopes. Recall that a {\it polytope} $P\subset\R^n$ is the convex hull of finitely many points in $\R^n$. Furthermore, $P$ is a {\it lattice} polytope if its vertices belong to the integer lattice $\Z^n\subset\R^n$.

Let $K$ be a convex body and $Q_2,\dots, Q_n$
be polytopes in $\R^n$. Given $u\in\mS^{n-1}$, let 
$V(Q^u_2,\dots, Q^u_{n})$ denote the $(n-1)$-dimensional mixed volume of 
$Q^u_2,\dots, Q^u_{n}$ translated to the orthogonal subspace~$u^\perp$. Then we have
\begin{equation}\label{E:induct}
V(K,Q_2\dots, Q_n)=\frac{1}{n}\sum_{u\in\mS^{n-1}}h_{K}(u)V(Q^u_{2},\dots,Q^u_{n}).
\end{equation}

Note that the above sum is finite, since there are only finitely many $u\in\mS^{n-1}$ for which 
$\{Q^u_{2},\dots,Q^u_{n}\}$ is essential. Namely, these $u$ are
among the outer unit normals to the facets of $Q_2+\dots+Q_n$.
 
 \begin{Rem}\label{R:induct-lattice} 
 There is a reformulation of \re{induct} that is more suitable for lattice polytopes. It is not hard to see that $n!\Vol(P)$
 is an integer for any lattice polytope. This implies that $n!V(P_1,\dots, P_n)$ is also an integer for any collection of lattice polytopes 
 $P_1,\dots, P_n$. 
 Recall that a vector $u\in\Z^n$ is {\it primitive} if the greatest common divisor of its components is 1. Given lattice polytopes $P, Q_2,\dots, Q_n$ we have
 \begin{equation}\label{E:induct-lattice}
n!V(P,Q_2\dots, Q_n)=\sum_{u\text{ primitive}}h_{P}(u)(n-1)!V(Q^u_{2},\dots,Q^u_{n}),
\end{equation}
where the $(n-1)$-dimensional mixed volume is normalized such that the volume of the parallelepiped spanned by a lattice basis for $u^\perp\cap\Z^n$ equals one. Note that the terms in the sum  are non-negative integers, which, as above, equal zero for all but finitely many primitive $u\in\Z^n$.
\end{Rem}

\subsection*{Cayley Polytopes and Combinatorial Cayley Trick}
Let $P_1,\ldots,P_k \subset  \R^n$ be convex polytopes. The associated {\it Cayley polytope}
$${\cC}(P_1,\ldots,P_k)$$
is the convex hull in $\R^n \times \R^k$ of the union of the polytopes $P_i {\times} \{e_i\}$ for $i=1,\ldots,k$, where $\{e_1,\ldots,e_k\}$ is the standard basis for $\R^k$. 

Let $(x,y)=(x_1,\ldots,x_n,y_1,\ldots,y_k)$ be coordinates on $\R^n \times \R^k$ and let $\pi_1: \R^n \times \R^k \rightarrow \R^n$
and $\pi_2: \R^n \times \R^k \rightarrow \R^k$ be the projections defined by $\pi_1(x,y)=x$ and $\pi_2(x,y)=y$, respectively. 
Note that $\pi_2({\cC}(P_1,\ldots,P_k))$ is the $(k-1)$-dimensional simplex $\D_{k-1}$ defined by $\sum_{i=1}^k y_i=1$ and $y_i\geq 0$ for $1\leq i\leq k$. 
Furthermore, for every $y\in \Delta_{k-1}$ we have
\begin{equation}\label{E:cayley}
\pi_2^{-1}(y) \cap {\cC}(P_1,\ldots,P_k)=(y_1P_1+\dots+y_kP_k)\times \{y\}.
\end{equation}
Note that when all $y_i>0$ the preimage $\pi_2^{-1}(y) \cap {\cC}(P_1,\ldots,P_k)$ has dimension equal to $\dim (P_1+\cdots+P_k)$. This implies that
\begin{equation}\label{E:dim-cayley}
\dim {\cC}(P_1,\ldots,P_k)=\dim (P_1+\cdots+P_k)+k-1.
\end{equation}
If $\dim P_i \geq 1$ for $i=1,\ldots,k$, then the Cayley polytope ${\cC}(P_1,\ldots,P_k) \subset \R^{n+k}$, as well as
the Minkowski sum $P_1+\cdots+P_k$, is called \emph{fully mixed}. The following result is an immediate consequence of \eqref{E:dim-cayley}.

\begin{Lemma}\label{basic}
Consider polytopes $P_1,\ldots,P_n \subset \R^n$. Then the following conditions are equivalent.
\begin{enumerate}
\item The Cayley polytope ${\cC}(P_1,\ldots,P_n)$ is a fully mixed $(2n-1)$-dimensional simplex.
\item $P_1, \ldots,P_n$ are segments with linearly independent directions.
\end{enumerate}
\end{Lemma}

\begin{Rem}\label{generalized} Let $P_1,\ldots,P_n$ be polytopes in $\R^n$. 
From \rt{essential} and Lemma~\ref{basic}, we have $V(P_1,\ldots,P_n)>0$ if and only if ${\cC}(P_1,\ldots,P_n)$
contains a fully mixed $(2n-1)$-dimensional simplex ${\cC}(E_1,\ldots,E_n)$.
\end{Rem}

Let $\tau_{\cC}$ be any polyhedral subdivision of ${\cC}(P_1,\ldots,P_k)$ with vertices in $ \cup_{i =1}^k P_i {\times}  \{e_i\}$.
Consider any full-dimensional polytope ${\cC}_{\sigma}$ of $\tau_{\cC}$. Then it intersects each hyperplane $y_i=1$ along a non-empty face $\sigma_i \times \{e_i\}\subset P_i\times\{e_i\}$ for $1\leq i\leq k$, and it follows that ${\cC}_{\sigma}={\cC}(\sigma_1,\ldots,\sigma_k)$. Therefore, $\tau_{\cC}$ consists of the set of all the polytopes ${\cC}(\sigma_1,\ldots,\sigma_k)$ together with their faces. Taking the image under $\pi_1$ of $\pi_2^{-1}\big(\frac{1}{k},\dots, \frac{1}{k}\big) \cap {\cC}(P_1,\ldots,P_n)$ we obtain, by \re{cayley}, the Minkowski sum $P_1+\cdots+P_k$ (up to dilatation by $\frac{1}{k}$) together
with a polyhedral subdivision by polytopes $\sigma_1+\cdots+\sigma_k$, where $\sigma_i \subset P_i$ for $1\leq i\leq k$. This defines a correspondence from the set of all polyhedral subdivisions of ${\cC}(P_1,\ldots,P_k)$ with vertices in $ \cup_{i =1}^k P_i {\times}  \{e_i\}$ to a set of polyhedral subdivisions of $P_1+\cdots+P_k$ which are called \emph{mixed}. Note that $\tau_{\cC}$ is uniquely determined by the corresponding mixed subdivision of $P_1+\cdots+P_k$. 
This one-to-one correspondence is commonly called the \emph{combinatorial Cayley trick} or simply the \emph{Cayley trick}, see \cite{HRS}, \cite{St} or \cite{BB}, for instance.

%

A mixed polyhedral subdivision of $P_1+\cdots+P_k$ is called {\it pure} if the corresponding subdivision of ${\cC}(P_1,\ldots,P_k)$ is a triangulation. Let $\sigma_1+\cdots+\sigma_k$ be a polytope in a pure mixed polyhedral subdivision of $P_1+\cdots+P_k$.
Then each $\sigma_i$ is a simplex since $\sigma_i \times \{e_i\}$ is a face of the simplex ${\cC}(\sigma_1,\ldots,\sigma_k)$.
If furthermore $k=n$ and $\dim (\sigma_1+\cdots+\sigma_n)=n$, then $\sigma_1+\cdots+\sigma_n$ is fully mixed if and only if ${\cC}(\sigma_1,\ldots,\sigma_n)$ is a fully mixed $(2n-1)$-dimensional simplex, equivalently,  $\sigma_1,\ldots,\sigma_n$ are segments with linearly independent directions (see Lemma \ref{basic}).
%
The following result is well-known, see \cite[Theorem 2.4]{HS} or \cite[Theorem 6.7]{USA}.

\begin{Lemma}\label{L:key}
For convex polytopes $P_1,\ldots,P_n$ in $\R^n$, the quantity $n! \, V(P_1,\dots, P_n)$ is equal to the sum of the Euclidean volumes of the fully mixed polytopes in any pure mixed polyhedral subdivision of $P_1+\cdots+P_n$.
\end{Lemma}

%

\section{First criterion}

In this section we present our first criterion for strict monotonicity of the mixed volume and its corollaries.

\begin{Def}
Let $K$ be a subset of a convex polytope $A$ and let $F\subset A$ be a facet.
We say {\it $K$ touches $F$} when the intersection $K\cap F$ is non-empty.
\end{Def}

We will often make use of the following proposition, which gives a criterion for strict monotonicity in a very special case, see \cite[page 282]{S}. 

\begin{Prop}\label{P:mix}
Let  $P_1, Q_1,\dots, Q_n$ be convex polytopes in $\R^n$ and $P_1\subseteq Q_1$.
Then $V(P_1,Q_2,\dots, Q_n)=V(Q_1,Q_2,\dots, Q_n)$ if and only if $P_1$ touches every face $Q_1^u$
for $u$ in the set
$$U=\{u\in \mS^{n-1}\ |\ \{Q_2^u,\dots, Q_n^u\}\ \text{\rm  is essential}\}.$$   
\end{Prop}

The above statement easily follows from \re{induct} and the observation $h_{P_1}(u)\leq h_{Q_1}(u)$ with equality if and only if $P_1$ touches $Q^u_1$. See \cite[Sec 5.1]{S}
for details.


Here is the first criterion for strict monotonicity.

\begin{Th}\label{T:main2} Let  $P_1,\dots, P_n$ and $Q_1,\dots, Q_n$ be convex polytopes in $\R^n$ such that $P_i\subseteq Q_i$ for every $i\in[n]$. Given $u\in\mS^{n-1}$ consider the set
$$T_u=\{ i\in[n] \ |\ P_i \text{ touches }Q^u_i \}.$$ 
Then $V(P_1,\dots, P_n) < V(Q_1,\dots, Q_n)$ if and only if there exists $u\in\mS^{n-1}$ such that the collection 
$\{Q^u_i\ |\ i\in T_u\}\cup \{Q_i\ |\ i\in [n]\setminus T_u\}$ is essential.
\end{Th}
%

\begin{pf} Assume that there exists $u\in\mS^{n-1}$ such that the collection 
$\{Q^u_i\ |\ i\in T_u\}\cup \{Q_i\ |\ i\in [n]\setminus T_u\}$ is essential.
Note that $T_u$ is a proper subset of $[n]$, otherwise $\{Q^u_i\ |\ i\in T_u\}$ is a collection of $n$ polytopes
contained in translates of an $(n-1)$-dimensional subspace, hence, cannot be essential. 
Without loss of generality we may assume
that $[n]\setminus T_u=\{1,\dots, k\}$ for some $k\geq 1$. In other words, we assume the collection
\begin{equation}\label{E:old}
\{{Q_1},\dots, {Q_k},Q^u_{k+1},\dots, Q^u_{n}\}
\end{equation}
is essential. Since $P_i$ does not touch $Q^u_i$ for $1\leq i\leq k$
there is a hyperplane $H=\{x\in\R^n\ |\ \la x,u\ra = h_{Q_i}(u)-\varepsilon\}$ which separates $P_i$ and $Q^u_i$.
Let $H_+$ be the half-space containing $P_i$. Then the truncated polytope $\tilde{Q}_i=Q_i\cap H_{+}$
satisfies $P_i\subseteq \tilde{Q}_i\subset Q_i$. We claim that, after a possible renumbering of the first $k$ of the $Q_i$, the collection 
\begin{equation}\label{E:new}
\{\tilde{Q^u_2},\dots, \tilde{Q^u_k},Q^u_{k+1},\dots, Q^u_{n}\}
\end{equation}
is essential. Indeed, since \re{old} is essential, by \rt{essential} there exist $n$ segments $E_i\subset Q_i$ with
linearly independent directions such that $E_i\subset Q^u_i$ for $k<i\leq n$. Replace the first $k$ of the segments with  
their projections onto $\tilde{Q^u_1},\dots, \tilde{Q^u_k}$. By \rl{linalg} below, after a possible renumbering of the first $k$ segments,
we obtain $n-1$ segments $E_2,\dots, E_n$ with linearly independent directions such that $E_i\subset \tilde{Q^u_i}$ for $2\leq i\leq k$ and $E_i\subset Q_i^u$ for $k<i\leq n$. By \rt{essential}, the collection \re{new} is essential.

Now, by \rp{mix} and since $P_1$ does not touch $Q_1^u$, we obtain
$$V(P_1,\tilde{Q_2},\dots, \tilde{Q_k},Q_{k+1},\dots, Q_{n})<V(Q_1,\tilde{Q_2},\dots, \tilde{Q_k},Q_{k+1},\dots, Q_{n}).$$
Finally, by monotonicity we have
$V(P_1,\dots, P_n)\leq V(P_1,\tilde{Q_2},\dots, \tilde{Q_k},Q_{k+1},\dots, Q_{n})$ and $V(Q_1,\tilde{Q_2},\dots, \tilde{Q_k},Q_{k+1},\dots, Q_{n})\leq V(Q_1,\dots, Q_n)$. Therefore,  $$V(P_1,\dots, P_n) < V(Q_1,\dots, Q_n).$$

Conversely, assume $V(P_1,\dots, P_n) < V(Q_1,\dots, Q_n)$. Then, by monotonicity, for some $1\leq k\leq n$ we have
$$V(P_1,\dots, P_{k-1},P_k,Q_{k+1},\dots, Q_{n})<V(P_1,\dots, P_{k-1},Q_k,Q_{k+1},\dots, Q_{n}).$$
By \rp{mix} there exists $u\in\mS^{n-1}$ such that $\{P^u_1,\dots, P^u_{k-1},Q^u_{k+1},\dots, Q^u_{n}\}$ is
essential and $k\not\in T_u$. By choosing a segment in $Q_k$ not parallel to the orthogonal hyperplane $u^\perp$ (which exists since $P_k\subset Q_k$, but $P_k$ does not touch $Q_k^u$)
we see that $$\{P^u_1,\dots, P^u_{k-1},Q_k,Q^u_{k+1},\dots, Q^u_{n}\}$$ is essential.
It remains to notice that $P^u_i\subseteq Q^u_i$ for $i\in T_u$ and, hence, the collection
$\{Q^u_i\ |\ i\in T_u\}\cup \{Q_i\ |\ i\in [n]\setminus T_u\}$ is essential as well.
\end{pf}

\begin{Lemma}\label{L:linalg} Let $\{v_1,\dots, v_k,v_{k+1},\dots, v_n\}$ be a basis for $\R^n$ where $v_{k+1},\dots,v_n$ belong to a hyperplane $H\subset\R^n$.
Let $\pi$ denote the orthogonal projection onto $H$. Then, after a possible renumbering of the first $k$ vectors, the set $\{\pi(v_2),\dots,\pi(v_k),v_{k+1},\dots,v_n\}$
is a basis for $H$. 
\end{Lemma}

\begin{pf} Clearly, the set $\{\pi(v_1),\pi(v_2),\dots,\pi(v_k),v_{k+1},\dots,v_n\}$ spans $H$. Starting with the linearly independent set $\{v_{k+1},\dots,v_n\}$ we can extend it to a basis for $H$ by appending $k-1$ vectors from $\{\pi(v_1),\dots,\pi(v_k)\}$.
\end{pf}

\begin{Rem}\label{R:full-dim}
Note that if $Q_1,\dots, Q_n$ are $n$-dimensional then $\{Q^u_i\ |\ i\in T_u\}\cup \{Q_i\ |\ i\in [n]\setminus T_u\}$ is essential if and only if $\{Q^u_i\ |\ i\in T_u\}$ is essential. (This can be readily seen from \re{essential}.) In this case we can simplify the criterion of \rt{main2} as follows:
$V(P_1,\dots, P_n) < V(Q_1,\dots, Q_n)$ if and only if there exists $u\in\mS^{n-1}$ such that the collection 
$\{Q^u_i\ |\ i\in T_u\}$ is essential.
\end{Rem}


\begin{Rem} After the initial submission of our paper to arxiv.org we were informed by Maurice Rojas that a 
similar criterion for rational polytopes appeared in his paper \cite[Corollary 9]{Rojas}. The proof of his criterion 
is algebraic and is based on the BKK bound. We are thankful to Maurice Rojas for pointing that out.
\end{Rem}

%

A particular instance of \rt{main2}, especially important for applications to polynomial systems, is the case when
$P_1,\dots, P_n$ are arbitrary polytopes and $Q_1,\dots, Q_n$ are equal to the same polytope $Q$. 
We will assume that $Q$ is $n$-dimensional, otherwise $\{P_1,\dots, P_n\}$
is not essential and, hence, both $V(P_1,\dots,P_n)$ and $V(Q,\dots,Q)=\Vol_n(Q)$ are zero.
Then the strict monotonicity has the following simple geometric interpretation.

\begin{Cor}\label{C:MV=V}
Let $P_1,\dots, P_n$  be polytopes in~$\R^n$ contained in an $n$-dimensional polytope $Q$.
Then $V(P_1,\dots, P_n)<\Vol_n(Q)$ if and only if 
there is a proper face of $Q$ of dimension $t$ which is touched by at most $t$ of the polytopes $P_1,\dots, P_n$.
\end{Cor}
\begin{pf}
By \rt{main2} and \rr{full-dim}, we have $V(P_1,\dots, P_n)<\Vol_n(Q)$ if and only if there exists $u\in\mS^{n-1}$ such that the collection $\{Q^u,\dots, Q^u\}$, where $Q^u$ is repeated $|T_u|$ times, is essential. The last condition is equivalent to
$\dim Q^u \geq |T_u|$. This precisely means that the face $Q^u$ is touched by at most $\dim Q^u$ of the polytopes.
\end{pf}
Note that in particular, if a vertex of $Q$ does not belong to any of the polytopes $P_1,\ldots,P_n$, then $V(P_1,\dots, P_n)<\Vol_n(Q)$.

\begin{Ex} \label{Ex:2-dimCor} Let $P_1,P_2$ be convex polytopes in $\R^2$ and $Q$ be the convex hull of their union.
Then \rc{MV=V} shows that  $V(P_1,P_2)<V(Q)$ if and only if either $P_1$ or $P_2$ does not touch some edge of 
$Q$.
\end{Ex} 

\begin{Rem} One can obtain a more direct proof of \rc{MV=V} by modifying the proof of 
Theorem 2.6 in \cite{SZ}. In the case when $Q$ equals the convex hull of $P_1\cup\dots\cup P_n$,
several sufficient conditions for $V(P_1,\dots, P_n)=\Vol_n(Q)$ were found
by Tianran Chen in \cite{Chen} whose preprint appeared on arxiv.org shortly after ours. Chen studies other
applications of \rc{MV=V} to sparse polynomial systems from the computational complexity point of view.
\end{Rem}

\begin{Ex} \label{Ex:first} Let $\cA$ be a finite set in $\R^n$ with $n$-dimensional convex hull $Q$
and choose a subset $\{a_1,\dots, a_n\}\subset\cA$. Define
$\cA_i=(\cA\setminus \{a_{1},\ldots, a_{n}\}) \cup \{a_{i}\}$ for $1\leq i\leq n$ and let $P_i$ be the convex hull of $\cA_i$.
Then \rc{MV=V} leads to $V(P_1,\ldots,P_n)=\Vol_n(Q)$. Indeed, for any $I \subset [n]$ of size $n-t$ and  any face $F\subset Q$ of dimension $t$ we have
$|\cup_{i \in I} \cA_i| = |\cA|-t$ and $|F \cap \cA| \geq t+1$. 
Therefore, the subsets $\cup_{i \in I} \cA_i$ and $F\cap \cA$ cannot be disjoint, i.e. every face $F$ of dimension $t$ is touched by the union of any $n-t$ of the $P_i$.
This shows that every face $F$ of dimension $t$ is touched by at least $t+1$ of the $P_i$ and we can apply \rc{MV=V}.
\end{Ex}
%

Going back to the statement of \rc{MV=V} it seems natural to ask: If there is a dimension $t$ face of $Q$ not touched by more than $n-t$ polytopes among $P_1,\dots,P_n$, can the inequality $V(P_1,\dots, P_n)<\Vol_n(Q)$ be improved?
The answer is clearly no in the class of all polytopes: For any $\varepsilon>0$, take $P_i=(1-\varepsilon)Q$ for $1\leq i\leq n$, then $$V(P_1,\dots, P_n)=(1-\varepsilon)^n\Vol_n(Q)$$
and no face of $Q$ is touched by any of the $P_i$. However, one might expect an improvement in the class of  { lattice polytopes}. 
We present such an improvement of \rc{MV=V} in \rp{stronger}.

Let $P\subset Q$ be lattice polytopes and $u\in\Z^n$ primitive. We can define the {\it lattice distance from $P$ to $Q$ in the direction of $u$} as the difference $h_{Q}(u)-h_{P}(u)$. Note that this is a non-negative integer. 

\begin{Prop}\label{P:stronger}
Let $P_1,\dots, P_n$  be lattice polytopes in~$\R^n$ contained in an $n$-dimensional lattice polytope $Q$. 
Suppose there exists a facet $Q^v\subset Q$, for some primitive $v\in\Z^n$, which is not touched by $P_1,\dots, P_{m}$, 
for some $m \geq 1$. Moreover, suppose that $\{P^v_1,\dots, P^v_{m-1}\}$ is essential.
Then $$n!\,\Vol Q-n!\,V(P_1,\dots, P_n)\geq l_1+\dots+l_m,$$
where $l_i\geq 1$ is the lattice distance from $P_i$ to $Q$ in the direction of $v$.
\end{Prop}
\begin{pf} By the essentiality of  $\{P^v_1,\dots, P^v_{m-1}\}$ and since $\dim Q^v=n-1$ it 
follows that for any $1\leq i\leq m$
the collection $\{P^v_1,\dots,P^v_{i-1},Q^v,\dots,Q^v\}$, where $Q^v$ is repeated $n-i$ times, is essential.
Applying \re{induct-lattice} from \rr{induct-lattice} we obtain

\begin{eqnarray}
& & n!V(P_1,\dots, P_{i-1},\underbrace{Q,\dots, Q}_{n-i+1})-n!V(P_1,\dots, P_{i},\underbrace{Q,\dots,Q}_{n-i})=\nonumber\\
& &\sum_{u\text{ primitive }}\left(h_{Q}(u)-h_{P_i}(u)\right)(n-1)!V(P^u_1,\dots,P^u_{i-1},Q^u,\dots,Q^u)\geq\nonumber\\
& &\left(h_{Q}(v)-h_{P_i}(v)\right)(n-1)!V(P^v_1,\dots,P^v_{i-1},Q^v,\dots,Q^v)\geq h_{Q}(v)-h_{P_i}(v)=l_i,\nonumber
\end{eqnarray}
where the last inequality follows from the fact that $\{P^v_1,\dots,P^v_{i-1},Q^v,\dots,Q^v\}$ is essential and, hence, $(n-1)!V(P^v_1,\dots,P^v_{i-1},Q^v,\dots,Q^v)$ is a positive integer. Summing up these inequalities for $1\leq i\leq m$ we produce
 $$n!\,\Vol Q-n!\,V(P_1,\dots, P_m,\underbrace{Q,\dots,Q}_{n-m})\geq l_1+\dots+l_m.$$
 Now the required inequality follows by the monotonicity of the mixed volume.
 
%
\end{pf}

Note that the case of $m=1$ and $l_1=1$ recovers an instance of \rc{MV=V} with $t=n-1$. In this case the condition that $\{P^u_1,\dots, P^u_{m-1}\}$ is essential is void. We remark that in general this condition cannot be removed. Indeed,
let $Q=\conv\{0,e_1,\dots,e_n\}$ be the standard $n$-simplex, $Q^u$ one of its facets, and $P_i\subset Q$ for $1\leq i\leq n$.
Then if $P_1,\dots, P_m$ equal the vertex of $Q$ not contained in $Q^u$, then 
$$0=n!V(P_1,\dots,P_n)=n!V_n(Q)-1,$$
regardless of $m$. It would be interesting to obtain a more general statement than \rp{stronger} which deals
with smaller dimensional faces, rather than facets.

\begin{Rem}
Given a polytope $Q$, \rc{MV=V} provides a characterization of collections $P_1,\dots,P_n$ such that
$P_i\subset Q$ for $1\leq i\leq n$ and $V(P_1,\dots,P_n)=\Vol_n(Q)$.  
Clearly, when $Q$ and the $P_i$ are lattice polytopes there are only finitely many such collections.
Describing them explicitly is a hard combinatorial problem, in general. In the case when $Q$ is the standard simplex we have
$\Vol_n(Q)=1/n!$. On the other hand, since  $n!V(P_1,\dots,P_n)$ is an integer when the $P_i$ are lattice polytopes, 
\rt{essential} implies the following: For lattice polytopes $P_1,\dots, P_n$ contained in the standard simplex $Q$ we have
 $V(P_1,\dots,P_n)=\Vol_n(Q)$ if and only if $\{P_1,\dots,P_n\}$ is essential. This is a particular case of a much deeper result 
by Esterov and Gusev, who described all collections of lattice polytopes $\{P_1,\dots, P_n\}$ satisfying $V(P_1,\dots,P_n)=1/n!$, see~\cite{EG}. 
\end{Rem}

%


\section{Second Criterion}\label{S:subd}
In this section, we obtain another criterion for the strict monotonicity property (\rt{main3}) based on mixed polyhedral 
subdivisions and the combinatorial Cayley trick. We first present a result
about faces of Cayley polytopes which will be useful. Consider convex polytopes $P_1,\ldots,P_k \subset \R^n$.
Recall that the Cayley polytope ${\cC}(P_1,\ldots, P_k)$ is the convex hull in $\R^{n+k}$ of the union of the polytopes
$P_i \times \{e_i\}$ for $1\leq i\leq k$, where $\{e_1,\ldots,e_k\}$ is the standard basis of $\R^k$ (see Section \ref{Pre}).
\begin{Lemma}\label{L:useful}
Let $(u,v)$ be a vector in $\R^{n} \times \R^k$ with $v=(v_1,\ldots,v_k)$.
Then $$h_{{\cC}(P_1,\ldots, P_k)}(u,v)=\max\{h_{P_i}(u)+v_i \ |\ i\in[k] \}.$$
Moreover, ${\cC}(P_1,\ldots, P_k)^{(u,v)}$ is the convex hull in $\R^{n}\times\R^{k}$
of the union of the polytopes $P_i^u\times \{e_i\}$ for $i$ in the set
$$I=\{i \in [k]\ |\ h_{{\cC}(P_1,\ldots, P_k)}(u,v)=h_{P_i}(u)+v_i\}.$$
\end{Lemma}
\begin{pf}
Let $(x,y) \in \R^n \times \R^k$ be a point of $ {\cC}(P_1,\ldots, P_k)$.
There exist $\alpha_i \in\R_{\geq 0}$ and points $x_i \in P_i$ for $i\in[k]$ such that $(x,y)=
\sum_{i=1}^k \alpha_i(x_i,e_i)$ and $\sum_{i=1}^k \alpha_i=1$.
Then $\la (u,v), (x,y) \ra=\sum_{i=1}^k \alpha_i (\la u, x_i \ra+v_i)$ is bounded above by
$\max\{h_{P_i}(u)+v_i \ |\  i\in[k] \}$. Moreover, this bound is attained if and only if
$\la u, x_i \ra+v_i=\max\{h_{P_i}(u)+v_i \ |\  i\in[k] \}$ for all $i$ such that $\alpha_i > 0$.
Since $\la u, x_i \ra+v_i \leq h_{P_i}(u)+v_i  \leq \max\{h_{P_i}(u)+v_i \ |\  i\in[k] \}$, the latter condition is equivalent to
$x_i \in P_i^u$ with $i \in I$.
\end{pf}

\begin{Rem}\label{hyperplane}
The polytope ${\cC}(P_1,\ldots, P_k)$ is contained in the hyperplane $\{(x,y) \in \R^n \times \R^k \ |\ \sum_{i=1}^k y_i=1\}$. Therefore,
$ {\cC}(P_1,\ldots, P_k)^{(u,v)}={\cC}(P_1,\ldots, P_k)$ for $u=0$ and $v= (\lambda,\ldots,\lambda)$ for any non-zero $\lambda \in \R$.
Moreover, if $\dim (P_1+\cdots+P_k)=n$, then ${\cC}(P_1,\ldots, P_k)^{(u,v)}={\cC}(P_1,\ldots, P_k)$ only if $u=0$
and $v= (\lambda,\ldots,\lambda)$  with $\lambda \in \R$, by~\eqref{E:dim-cayley}.
\end{Rem}
Before stating the main result of this section, we need a technical lemma about regular polyhedral subdivisions
(see, for example, \cite{HRS} for a reference on this topic).
Let $\mathcal A$ be a finite set in $\R^n$ and let $A$ denote the convex hull of $\mathcal A$. A polyhedral subdivision $\tau$ of
$A$ with vertices in $\mathcal A$ is called \emph{regular} if there exists a map $h:{\mathcal A} \rightarrow \R$ such that $\tau$ is obtained by projecting the lower faces of the convex-hull $\hat{A}$ of $\{(a,h(a)) \ |\ a \in {\mathcal A}\}$ via the projection $\R^{n+1} \rightarrow \R^n$ forgetting the last coordinates. Here a lower face of $\hat{A}$  is a face of a facet of $\hat{A}$ with an inward normal vector with positive last coordinate. Intuitively, this is a face that can be seen from below in the direction of $e_{n+1}$.The union of the lower faces of $\hat{A}$ is the graph of a convex piecewise linear map $\hat{h}: A \rightarrow \R$ whose domains of linearity are the polytopes of $\tau$. We say that $h$ (respectively $\hat{h}$) \emph{certifies the regularity} of $\tau$ and that $\tau$ is induced by $h$ (respectively $\hat{h}$). Note that for $h$ generic enough any $n+2$ points of $\{(a,h(a)) \ |\ a \in {\mathcal A}\}$ are affinely independent (i.e. do not lie on a hyperplane), hence the induced subdivision $\tau$ is a triangulation.
%
%

\begin{Lemma}\label{regular}
Let ${\mathcal A}_1, {\mathcal A}_2$ be finite subsets of $\R^n$ and $A_1, A_2$ their convex hulls.
\begin{enumerate}
\item If ${\mathcal A}_1 \subset {\mathcal A}_2$ and if $\tau_1$ is a regular polyhedral subdivision (respectively a regular triangulation)
of  $A_1$ with vertices in
${\mathcal A}_1$,  then there exists a regular polyhedral subdivision (respectively a regular triangulation) $\tau_2$ of  $A_2$ with vertices in
${\mathcal A}_2$ such that $\tau_1 \subset \tau_2$.
\item Assume that the relative interiors of $A_1$ and $A_2$ do not intersect. Let ${\mathcal A}={\mathcal A}_1 \cup {\mathcal A}_2$ and let $A$ denote the convex hull of ${\mathcal A}$.
If $\tau_1$ is a regular polyhedral subdivision (respectively a regular triangulation) of  $A_1$ with vertices in ${\mathcal A}_1$ and $\tau_2$ is a regular polyhedral subdivision (respectively a regular triangulation) of  $A_2$ with vertices in ${\mathcal A}_2$, then there exists a regular polyhedral subdivision (respectively a regular triangulation) $\tau$ of  $A$ with vertices in ${\mathcal A}$
such that $\tau_1 \cup \tau_2 \subset \tau$.
\end{enumerate}
\end{Lemma}
\begin{pf}
(1) Consider a map $H: {\mathcal A}_2 \rightarrow \R$ which vanishes on ${\mathcal A}_1$ and
takes positive values on ${\mathcal A}_2 \setminus {\mathcal A}_1$. Then, this map certifies the regularity of a polyhedral subdivision $\tilde{\tau}_2$ of $A_2$ with vertices in ${\mathcal A}_2$ and which contains $A_1$. Now consider a regular polyhedral subdivision
$\tau_1$ of $A_1$ with vertices in ${\mathcal A}_1$ whose regularity is certified by
$h_1:{\mathcal A}_1 \rightarrow \R$. Then, for $\epsilon>0$ small enough, the function $h_2:{\mathcal A}_2 \rightarrow \R$
defined by $h_2(a)=H(a)+\epsilon h_1(a)$ if $a \in {\mathcal A}_1$ and by $h_2(a)=H(a)$ otherwise certifies a regular polyhedral subdivision $\tau_2$ of $A_2$ with vertices in ${\mathcal A}_2$ such that $\tau_1 \subset \tau_2$.
Finally, if $\tau_1$ is a triangulation and if the values of $H$ on ${\mathcal A}_2 \setminus {\mathcal A}_1$ are generic enough, then $\tau_2$ is a triangulation.

(2) Consider a regular polyhedral subdivision (respectively a regular triangulation) $\tau_i$ of $A_i$ with vertices in ${\mathcal A}_i$ for $i=1,2$. 
Let $h_i : {\mathcal A}_i \rightarrow \R$ be a function certifying the regularity of $\tau_i$.
Since the relative interiors of the convex sets $A_1$ and $A_2$ do not intersect, there is a hyperplane which separates $A_1$ and $A_2$.
This means that there exist $u \in \R^n$ and $c \in \R$ such that $A_1 \subset B_-=\{x \in \R^n \ | \ \langle u , x \rangle  \leq c \}$, $A_2
\subset B_+=\{x \in \R^n \ | \ \langle u , x \rangle  \geq c \}$ and $A_1,A_2$ are not both contained in the separating hyperplane $B_+ \cap B_-$. Let $H: {\mathcal A} \rightarrow \R$ be a piecewise linear map which vanishes on $B_-$ and positive on $B_+ \setminus B_-$. For $\epsilon>0$ consider $h: {\mathcal A} \rightarrow \R$ defined by $h(a)=H(a)+\epsilon h_i(a)$ if $a \in {\mathcal A}_i$. If $\epsilon$ is small enough and generic, the map $h$ certifies the regularity of a polyhedral subdivision (respectively a regular triangulation) $\tau$ of $A$ such that $\tau_1 \cup \tau_2 \subset \tau$.
\end{pf}

Recall that for convex polytopes $Q_1,\ldots,Q_n$ in $\R^n$, we have $V(Q_1,\ldots,Q_n) >0$ if and only if the collection $\{Q_1,\ldots,Q_n\}$ is essential, which is equivalent to the existence of a fully mixed $(2n-1)$-dimensional simplex ${\cC}(E_1,\ldots,E_n)$ contained in ${\cC}(Q_1,\ldots,Q_n)$,
see \rt{essential} and Remark \ref{generalized}. We now describe a generalization of these equivalences. Consider convex
polytopes $P_i \subseteq Q_i \subset \R^n$ for $1\leq i\leq n$. For any non-zero vector $u \in \R^{n}$ define convex polytopes

$$B_{i,u}=\{x \in Q_i \ | \ \langle u,x \rangle \geq
h_{P_i}(u) \big\} , \ 1\leq i\leq n.$$
Intuitively, $B_{i,u}$ is the part of $Q_i$ lying on top of $P_i$ if one looks in the direction of the vector $u$.
%
%
%

\begin{Th} \label{T:main3}
Let  $P_1,\dots, P_n$ and $Q_1,\dots, Q_n$ be convex polytopes in $\R^n$ such that $P_i\subseteq Q_i$ for every $i\in[n]$.
The following conditions are equivalent:
\smallskip

\begin{enumerate}
\item $V(P_1,\ldots,P_n) < V(Q_1,\ldots,Q_n)$,

\item there exists a fully mixed $(2n-1)$-dimensional simplex ${\cC}(E_1,\ldots,E_n)$ contained in
the relative closure of ${\cC}(Q_1,\ldots,Q_n) \setminus {\cC}(P_1,\ldots,P_n)$,

\item there exists a non-zero vector $u \in \R^{n}$ such that the collection $\{B_{1,u},\ldots,B_{n,u}\}$ is essential.

%
%
%
%
%
%
%
%
\end{enumerate}
\end{Th}

%
%
%
%
%
%
%
\begin{pf}

First we note that if $\dim(Q_1+\cdots+Q_n) < n$, then none of the conditions (1), (2) and (3) holds.
Indeed, if  $\dim(Q_1+\cdots+Q_n) < n$, then $V(P_1,\ldots,P_n)=V(Q_1,\ldots,Q_n)=0$. Moreover,
$\dim{\cC}(Q_1,\ldots,Q_n) < 2n-1$ by \eqref{E:dim-cayley} and, thus, ${\cC}(Q_1,\ldots,Q_n)$
cannot contain a $(2n-1)$-dimensional simplex. Finally, (3) does not hold since otherwise
$Q_1+\cdots+Q_n$ would contain a fully mixed polytope which has dimension $n$. When $(P_1,\ldots,P_n)=(Q_1,\ldots,Q_n)$,
we also conclude that  none of the conditions (1), (2) and (3) holds for obvious reasons.

Assume now that $(P_1,\ldots,P_n) \neq (Q_1,\ldots,Q_n)$ and $\dim(Q_1+\cdots+Q_n)=n$. Write ${\bf (0,1)}$ for the vector
$(u,v)$ with $u=(0,\ldots,0) \in \R^n$ and $v=(1,\ldots,1) \in \R^n$.
Then ${\cC}(Q_1,\ldots,Q_n)$ has dimension $2n-1$ and its affine span
is a hyperplane orthogonal to ${\bf (0,1)}$, see Remark \ref{hyperplane}. Consider a fully mixed simplex ${\cC}_{E} ={\cC}(E_1,\ldots,E_n) \subset {\cC}(Q_1,\ldots,Q_n)$. Here $E_1,\ldots,E_n$ are segments with linearly independent directions contained in $Q_1,\ldots,Q_n$, respectively, and ${\cC}_{E}$ is the convex hull of the union $\cup_{i=1}^n E_i \times \{e_i\}$, by Lemma \ref{basic}.
Since ${\cC}(P_1,\ldots,P_n)$ and ${\cC}_{E}$ are convex sets contained in a hyperplane orthogonal to ${\bf (0,1)}$, the simplex
${\cC}_{E}$ is contained in the relative closure of ${\cC}(Q_1,\ldots,Q_n) \setminus {\cC}(P_1,\ldots,P_n)$ if and only
if there is a vector $(u,v) \in \R^n \times \R^n$ such that $(u,v)$ and ${\bf (0,1)}$ are not collinear and
%
%
%
%
\begin{equation}\label{E:contained}
{\cC}_{E} \subset \{(x,y) \in \R^n \times \R^n \ | \ \langle u,x \rangle + \langle v,y \rangle \geq h_{{\cC}(P_1,\ldots,P_n)}(u,v)\}.
\end{equation}
Note that the hyperplane defined by $\langle u,x \rangle + \langle v,y \rangle = h_{{\cC}(P_1,\ldots,P_n)}(u,v)$ is
a supporting hyperplane of ${\cC}(P_1,\ldots,P_n)$ and \eqref{E:contained} means that this hyperplane separates ${\cC}(P_1,\ldots,P_n)$ and ${\cC}_{E}$. Recall that $\dim {\cC}_{E} =2n-1$ and thus ${\cC}_{E}$ cannot be contained in this supporting hyperplane since otherwise it would be contained in the intersection of two distinct hyperplanes. Since $C_E$ is the convex-hull of the union of the polytopes $E_i \times \{e_i\}$ for $1\leq i\leq n$, we see that
\eqref{E:contained} is equivalent to
\begin{equation}\label{E:inter}
E_i \subset \{x \in \R^{n} \ | \ \langle u,x \rangle  \geq h_{{\cC}(P_1,\ldots,P_n)}(u,v)-v_i,\ 1\leq i\leq n\}.
\end{equation}
By Lemma \ref{L:useful}, we have $h_{{\cC}(P_1,\ldots, P_n)}(u,v)=\max\{h_{P_i}(u)+v_i \ |\ 1\leq i\leq n \}$. For
$u=0$ equation \eqref{E:inter} implies $v_1=\cdots=v_n$, which contradicts the fact that $(u,v)$ and ${\bf (0,1)}$ are not collinear (note that $h_{P_i}(u)=0$ when $u=0$). Therefore, if \eqref{E:inter} holds, then $u \neq 0$. Moreover, we get $E_i \subset B_{i,u}$ for $1\leq i\leq n$ since $h_{{\cC}(P_1,\ldots,P_n)}(u,v)-v_i \geq h_{P_i}(u)$. Consequently, \eqref{E:inter} implies that $\{B_{1,u},\ldots,B_{n,u}\}$ is essential, as the $E_i$ have linearly independent directions (see \rt{essential}.) We have proved the implication (2) $\Rightarrow $ (3).

Assume (3) holds and let $u \in \R^n$ be a non-zero vector such that $\{B_{1,u},\ldots,B_{n,u}\}$ is essential.
Let $E_i \subset B_{i,u}$, for $1\leq i\leq n$, be segments with linearly independent directions and choose $v=(v_1,\ldots,v_n) \in \R^n$
such that $h_{P_1}(u)+v_1=\cdots=h_{P_n}(u)+v_n$. Then, by \rl{useful}, $h_{{\cC}(P_1,\ldots,P_n)}(u,v)=h_{P_i}(u)+v_i$ for $1\leq i\leq n$ and
\eqref{E:inter} follows from $E_i \subset B_{i,u}$, $1\leq i\leq n$. Since \eqref{E:inter} is equivalent to (2), we conclude that (3) $\Rightarrow $ (2).
%
%
%
%
%
%
%
%
%
%
%


Now assume that (2) holds. Then, there exists a fully mixed $(2n-1)$-dimensional simplex $C_E$ satisfying \eqref{E:contained}.
It follows then from Lemma \ref{regular} that there exists a triangulation of ${\cC}(Q_1,\ldots,Q_n)$
with vertices in $ \cup_{i =1}^n Q_i { \times} \{e_i\}$ which contains ${\cC}_{E}$ and restricts to a triangulation of ${\cC}(P_1,\ldots,P_n)$.
Indeed, we may apply part (2) of Lemma~\ref{regular} to the set of vertices ${\mathcal A}_1$ of $C_E$ and to the set of vertices ${\mathcal A}_2$ of ${\cC}(P_1,\ldots,P_n)$. Taking for $\tau_1$ the trivial triangulation of $C_E$ and for $\tau_2$ any regular triangulation induced by a generic function $h_2$, we get the existence of a regular triangulation of the convex hull of ${\mathcal A}={\mathcal A}_1 \cup {\mathcal A}_2$ which contains $C_E$ and $\tau_2$. Applying part (1) of Lemma \ref{regular} to ${\mathcal A} \subset {\mathcal A} \cup V$ where $V$ is the set of vertices of $C(Q_1,\ldots,Q_n)$ gives a regular triangulation of ${\cC}(Q_1,\ldots,Q_n)$ as required.
By the combinatorial Cayley trick, this corresponds to a pure mixed subdivision $\tau_Q$ of $Q_1+\cdots+Q_n$ restricting to a pure mixed subdivision $\tau_P$ of $P_1+\cdots+P_n$ and a fully mixed polytope $E$ contained in
$\tau_Q \setminus \tau_P$. Therefore, $V(P_1,\dots,P_n) < V(Q_1,\ldots,Q_n)$ by Lemma \ref{L:key}.
This proves the implication $(2) \Rightarrow (1)$.

Assume now that  $V(P_1,\dots,P_n) < V(Q_1,\ldots,Q_n)$. Then, as in the proof of Theorem \ref{T:main2}, there exist $u\in\mS^{n-1}$ and $1\leq k\leq n$
such that $\{P^u_1,\dots, P^u_{k-1},Q^u_{k+1},\dots, Q^u_{n}\}$ is essential and $P_k$ does not touch $Q_k^u$. By choosing a segment in $B_{k,u}$ not parallel to the hyperplane $u^\perp$, which exists since $P_k$ does not touch $Q_k^u$, we conclude that 
$\{P^u_1,\dots, P^u_{k-1},B_{k,u},Q^u_{k+1},\dots, Q^u_{n}\}$ is essential as well. It remains to note that $P^u_i$ and 
$Q^u_{i}$ are contained in $B_{i,u}$ for $i\neq k$, $1\leq i\leq n$. We have proved the implication $(1) \Rightarrow (2)$.
\end{pf}

\begin{Rem}
Note that if $P_i$ touches $Q^u_i$ then $B_{i,u}=Q^u_i$ and if $P_i$ does not touch $Q^u_i$ then 
$\dim B_{i,u}=\dim Q_i$. Therefore, the condition (3) in the above theorem is equivalent to the condition 
in \rt{main2}.
\end{Rem}


\section{Polynomial systems}\label{S:pol}

Consider a finite set $\cA=\{a_1,\ldots,a_{\ell}\} \subset \Z^n$ where $\ell=|\cA|$. Let $(a_{1j},\ldots,a_{nj})$ be the coordinates of $a_j$ for $1\leq j\leq \ell$.
Consider a Laurent polynomial system with coefficients in an {algebraically closed field} $\K$
\begin{equation}\label{E:system}
f_1(x)=\cdots=f_n(x)=0,
\end{equation}
where $f_i(x)=\sum_{j=1}^{\ell} c_{ij}x^{a_j}$ and $x^{a_j}$ stands for the monomial $x_1^{a_{1j}} \cdots x_n^{a_{nj}}$. We 
assume that no polynomial $f_i$ is the zero polynomial. Call $\cA_i=\{a_j \in \cA \ | \ c_{ij} \neq 0\}$ the {\it individual support} of $f_i$. We may assume that
for any $1\leq j\leq \ell$ there exists $i$ such that $c_{ij} \neq 0$. Then $\cA=\cup_{i=1}^n \cA_i$ is called the {\it total support} of the system \eqref{E:system}. The {\it Newton polytope} $P_i$ of $f_i$ is the convex hull of $\cA_i$ and the {\it Newton polytope} $Q$ of the system \eqref{E:system} 
is the convex hull of $\cA$.

The matrices
$$C=(c_{ij}) \in \K^{n \times {\ell}} \; \mbox{and} \; A=(a_{ij}) \in \Z^{n \times {\ell}}$$
are the \emph{coefficient} and \emph{exponent} matrices of the system, respectively.
%

Choose $u\in\mS^{n-1}$ and let $\cA_i^u=P^u_i\cap\cA_i$. Then the {\it restricted system} corresponding to $u$
is the system
$$
f^u_1(x)=\cdots=f^u_n(x)=0,
$$
where $f^u_i(x)=\sum_{j=1}^{\ell} c^u_{ij}x^{a_j}$ with $c^u_{ij}=c_{ij}$ if $a_j\in\cA^u_i$ and $c^u_{ij}=0$, otherwise.
Finally, a system \re{system} is called {\it non-degenerate} if for every $u\in\mS^{n-1}$ the corresponding restricted system is inconsistent.

The relation between mixed volumes and polynomial systems originates in the following fundamental result, known as the BKK bound, discovered by Bernstein, Kushnirenko, and Khovanskii, see  \cite{Be, Kho, Kush}.

\begin{Th}\label{T:BKK}
The system \eqref{E:system} has at most $n!\,V(P_1,\ldots,P_n)$ isolated solutions in $(\K^*)^n$ counted with multiplicity. Moreover, it has precisely
$n!\,V(P_1,\ldots,P_n)$ solutions in $(\K^*)^n$ counted with multiplicity if and only if it is non-degenerate.
\end{Th}

\begin{Rem}
Systems with fixed individual supports and generic coefficients are non-degenerate. Moreover, the non-degeneracy condition is not needed if one passes to the toric compactification $X_P$ associated with the polytope $P=P_1+\dots+P_n$. Namely, a system has at most $n!\,V(P_1,\ldots,P_n)$ isolated solutions in $X_P$ counted with multiplicity, and if it has a finite number of solutions in $X_P$ then this number equals $n! \,V(P_1,\ldots,P_n)$ counted with multiplicity.
\end{Rem}

There are two operations on \eqref{E:system} that preserve its number of solutions in the torus $(\K^*)^n$: Left multiplication of $C$ by an
element of $\mbox{GL}_n(\K)$ and left multiplication of the {\it augmented exponent matrix}
$$
\bar{A}=\left(
\begin{array}{ccc}
1 & \cdots & 1 \\
& A & 
\end{array}
\right)
\in \Z^{(n+1) \times \ell}$$
by a matrix in $\mbox{GL}_{n+1}(\Z)$ whose first row is $(1,0,\dots, 0)$.
The first operation obviously produces an equivalent system. The second operation amounts
to applying an invertible affine transformation with integer coefficients  $a \mapsto b+\ell(a)$ on the total support $\cA
\subset \Z^n$. Here $b\in\Z^n$ is a translation vector and $\ell:\R^n \rightarrow \R^n$ is a linear map whose matrix
with respect to the standard basis belongs to $\mbox{GL}_{n}(\Z)$. A basic result of toric geometry says that there
is a monomial change of coordinates $x \mapsto y$ of the torus $(\K^*)^n$ so that $x^a=y^{\ell(a)}$ for any $a \in (\K^*)^n$.
Moreover, translating $\cA$ by $b$ amounts to multiplying each equation of \eqref{E:system} by the monomial $x^b$.
Thus starting from \eqref{E:system}, we get a system with same coefficient matrix, total support $b+\ell(\cA)$, and the same number
of solutions in $(\K^*)^n$.

\begin{Rem}\label{R:left}
Consider a non-degenerate system with coefficient matrix $C$.
While left multiplication of $C$ by an invertible matrix does not preserve the individual supports and Newton polytopes in general,
it preserves the total support of the system. Indeed, since $\cA$ is the total support of the system, no column of $C$ is zero, and thus no column of $C$ can become zero after left multiplication by an invertible matrix.
\end{Rem}
%
%

\begin{Ex}
Assume that \eqref{E:system} has precisely $n! \, \Vol_n(Q)$ solutions in $(\K^*)^n$ counted with multiplicity and $C$ has a non-zero maximal minor. Up to renumbering,  we may assume that this minor is given by the first $n$ columns of $C$. Left multiplication of \eqref{E:system} by the inverse of the corresponding submatrix of $C$ gives an equivalent system with invidual supports $\cA_i' \subset \cA_i''=(\cA\setminus \{a_{1},\ldots, a_{n}\}) \cup \{a_{i}\}$ for $1\leq i\leq n$. Thus, this new system has precisely $n!\, \Vol_n(Q)$ solutions in $(\K^*)^n$ counted with multiplicity. By Theorem \ref{T:BKK} this number of solutions is at most $n! \, V(P'_1,\ldots,P'_n)$, where $P_i'$ is the convex hull of $\cA'_i$. On the other hand, by monotonicity of the mixed volume we have $V(P'_1,\ldots,P'_n) \leq V(P_1'',\ldots,P_n'') \leq \Vol_n(Q)$. We conclude that $V(P'_1,\ldots,P'_n) = V(P_1'',\ldots,P_n'')= \Vol_n(Q)$. The second equality is also a consequence of Corollary \ref{C:MV=V}, see Example \ref{Ex:first}.
\end{Ex}

\begin{Th}\label{T:Ber}
Assume $\dim Q=n$.
If a system \eqref{E:system} has $n!\Vol_n(Q)$ isolated solutions in $(\K^*)^n$ counted with multiplicity, then for any proper face $F$ of $Q$ the submatrix $C_{\cF} \in \K^{n \times |\cF|}$ of $C$ with columns indexed by $\cF=\{j \in [{\ell}] \, , \, a_j \in F \cap \cA\}$
satisfies
$$
\text{\rm rank} \,  C_{\cF} \geq \dim F +1,
$$
or equivalently,
\begin{equation}\label{E:rank}
\text{\rm rank} \, C_{\cF} \geq \text{\rm rank} \, \bar{A}_{\cF}.
\end{equation}

Conversely, if \eqref{E:rank} is satisfied for all proper faces $F$ of $Q$ and if the system \eqref{E:system} is non-degenerate, then it has precisely $n!\Vol_n(Q)$ isolated solutions in $(\K^*)^n$ counted with multiplicity.
\end{Th}
\begin{pf}
First, note that for any proper face $F$ of $Q$ we have $\text{\rm rank} \, \bar{A}_{\cF}= \dim F+1$.
Consider a proper face $F$ of $Q$ of codimension $s \geq 1$ and assume that
$\text{\rm rank} \, C_{\cF} \leq \dim F=n-s$. Then there exist an invertible matrix
$L$ and $I \subset [n]$ of size $|I|=s$ such that the submatrix of $C'=LC$ with rows indexed by $I$ and columns indexed by ${\cF}$ is the zero matrix. The matrix $C'$ is the coefficient matrix of an equivalent system with the same total support, see Remark \ref{R:left}. Denote by $P_1',\ldots,P_n'$ the individual Newton polytopes of this equivalent system.
Then the polytopes $P_i'$ for $i \in I$ do not touch the face $F$ of $Q$, as $I$ corresponds to the zero submatrix of $C'$. Since $\dim F=n-s$ and $|I|=s$, it follows then from \rc{MV=V} that $V(P_1',\ldots,P_n') < \Vol_n(Q)$.
Theorem \ref{T:BKK} applied to the system with coefficient matrix $C'$ gives that it has at most $n! V(P_1',\ldots,P_n') < n! \Vol_n(Q)$ isolated solutions in $(\K^*)^n$ counted with multiplicity.
The same conclusion holds for the equivalent system \eqref{E:system}.
Therefore, if \eqref{E:system} has $n! \Vol_n(Q)$ isolated solutions in $(\K^*)^n$ counted with multiplicity, then \eqref{E:rank} is satisfied for all proper faces $F$ of $Q$.

Conversely, assume that \eqref{E:system} is non-degenerate and that \eqref{E:rank} is satisfied for all proper faces $F$ of $Q$. Then \eqref{E:system}  has $n! \ V(P_1,\ldots,P_n)$ isolated solutions in $(\K^*)^n$ counted with multiplicity by Theorem \ref{T:BKK}.
Suppose that $V(P_1,\ldots,P_n) < \Vol_n(Q)$. Then by \rc{MV=V}  there exists a proper face $F$ of $Q$ of codimension $s \geq 1$ and $I \subset [n]$ of size $|I|=s$ such that the polytopes
$P_i$ for $i \in I$ do not touch $F$. But then $\text{\rm rank} \, C_{\cF} \leq n-s=\dim F$, which gives a contradiction. Thus $V(P_1,\ldots,P_n) =\Vol_n(Q)$ and \eqref{E:rank}  has $n! \Vol_n(Q)$  isolated solutions in $(\K^*)^n$ counted with multiplicity.
\end{pf}

As an immediate consequence of Theorem \ref{T:Ber} from which we keep the notation, we get  the following corollary.

\begin{Cor} Consider any Laurent polynomial system \eqref{E:system} with $\dim Q=n$.
If there exists a proper face $F$ of $Q$ such that $\text{\rm rank} \, C_{\cF} < \text{\rm rank} \, \bar{A}_{\cF}$, then the system has either infinitely many solutions or strictly less than $n!\Vol_n(Q)$ solutions in $(\K^*)^n$ counted with multiplicity.
\end {Cor}
\begin{pf} 
Assume the existence of a proper face $F$ of $Q$ such that $
\text{\rm rank} \, C_{\cF} < \text{\rm rank} \, \bar{A}_{\cF}$. If \eqref{E:system} has precisely $n!\Vol_n(Q)$ solutions in $(\K^*)^n$ counted with multiplicity, then it is non-degenerate by \rt{BKK} and thus $\text{\rm rank} \, C_{\cF} \geq  \text{\rm rank} \, \bar{A}_{\cF}$ by \rt{Ber}, a contradiction.
\end{pf}
%
%

%
%


A very nice consequence of Theorem \ref{T:Ber} is the following result, which can be considered as a generalization of
Cramer's rule to polynomial systems.

\begin{Cor}\label{C:nice}
Assume that $\dim Q=n$ and that no maximal minor of $C$ vanishes. Then
the system \eqref{E:system} has the maximal number of $n!\,\Vol_n(Q)$ isolated solutions in $(\K^*)^n$ counted with multiplicity.
\end {Cor}
%
%
%
\begin{pf}
Note that $\ell \geq n+1$ since $\dim Q=n$ (recall that $\ell$ is the number of columns of $C$). Thus a maximal minor of $C$ has size $n$
and the fact that no maximal minor of $C$ vanishes implies that for any $J \subset [\ell]$ the submatrix of $C$ with rows indexed by $[n]$ and columns indexed by $J$ has maximal rank.
This rank is equal to $n$ if $|J| \geq n$ or to $|J|$ if $|J|< n$. Since $|\cF| \geq \dim F+1=\text{\rm rank} \, \bar{A}_{\cF}$ for any face $F$ of $Q$,
we get that $\text{\rm rank} \, C_{\cF} \geq \text{\rm rank} \, \bar{A}_{\cF}$ for any proper face $F$ of $Q$. Moreover, no restricted system is consistent for otherwise this would give a non-zero vector in the kernel of the corresponding submatrix of $C$. Thus \eqref{E:system} is non-degenerate and the result follows from Theorem \ref{T:Ber}. 
\end{pf}

When the polytope $Q=\conv\{0,e_1,\ldots,e_n\}$ is the standard simplex, the system \eqref{E:system} is linear
and it has precisely $n!\,\Vol_n(Q)=1$
solution in $(\K^*)^n$ if and only if no maximal minor of $C \in \K^{n \times (n+1)}$ vanishes, in accordance with
Cramer's rule for linear systems.

\section{Examples}

We conclude with a few examples illustrating the results of the previous section.

\begin{Ex} Let $\cA_1=\{(0,0),(1,2),(2,1)\}$ and $\cA_2=\{(2,0),(0,1),(1,2)\}$ be individual supports, and
$\cA=\cA_1\cup\cA_2$ the total support of a system. The Newton polytopes $P_1=\conv\cA_1$, $P_2=\conv\cA_2$,
and $Q=\conv\cA$ are depicted in \rf{ex1}, where the vertices of $P_1$ and $P_2$ are labeled by $\{1,2,3\}$ and $\{4,5,2\}$, respectively.
\begin{figure}[h]
\includegraphics[scale=1.8]{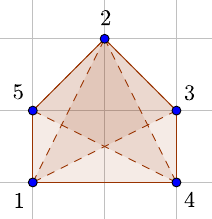}
\caption{The mixed volume of the two triangles equals the volume of the pentagon.}
\label{F:ex1}
\end{figure}
We use the labeling in \rf{ex1} to order the columns of the augmented matrix
$$\bar A=\left(\begin{matrix}
1 & 1 & 1 & 1 & 1\\
0 & 1 & 2 & 2 & 0\\
0 & 2 & 1 & 0 & 1 
\end{matrix}\right).$$
A general system with these supports has the following coefficient matrix
$$C=\left(\begin{matrix}
c_{11} & c_{12} & c_{13} & 0 & 0\\
0 & c_{22} & 0 & c_{24} & c_{25}
\end{matrix}\right),$$
where $c_{ij}\in\K$ are non-zero.
One can see that each edge of $Q$ is touched by both $P_1$ and $P_2$, hence, $V(P_1,P_2)=\Vol_2(Q)$,
see \rex{2-dimCor}. Also one can check that the rank conditions $
\text{\rm rank} \, C_{\cF} \geq \text{\rm rank} \, \bar{A}_{\cF}$ are satisfied for every face $F$ of $Q$. (In fact, both
ranks equal 2 when $F$ is an edge and 1 when $F$ is a vertex.)
\end{Ex}

\begin{Ex} Now we modify the previous example slightly, keeping $\cA_1$ the same and changing one of the
points in $\cA_2$, so $\cA_2=\{(2,0),(0,1),(1,1)\}$, see \rf{ex2}. 
\begin{figure}[h]
\includegraphics[scale=1.8]{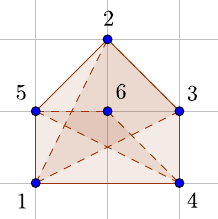}
\caption{The mixed volume of the two triangles is less than the volume of the pentagon.}
\label{F:ex2}
\end{figure}
The augmented exponent matrix and the coefficient matrix are as follows.
$$\bar A=\left(\begin{matrix}
1 & 1 & 1 & 1 & 1 & 1\\
0 & 1 & 2 & 2 & 0 & 1\\
0 & 2 & 1 & 0 & 1 & 1 \\
\end{matrix}\right),\quad
C=\left(\begin{matrix}
c_{11} & c_{12} & c_{13} & 0 & 0 & 0\\
0 & 0 & 0 & c_{24} & c_{25} & c_{26}
\end{matrix}\right).$$
This time the edge of $Q$ labeled by $\cF=\{2,3\}$ is not touched by $P_2$ and, hence, $V(P_1,P_2)<\Vol_2(Q)$.
Also, the rank condition for $\cF=\{2,3\}$ fails: $\text{\rm rank} \, C_{\cF}=1$ and $\text{\rm rank} \, \bar{A}_{\cF}=2$.
\end{Ex}

\begin{Ex}\label{Ex:prism}
Consider a system defined by the following 
augmented exponent matrix and coefficient  matrix
$$\bar A=\left(\begin{matrix}
\ 1 & 1 & 1 & 1 & 1 & 1\ \\
\ 0 & 1 & 1 & 0 & 0 & 0 \ \\
\ 0 & 0 & 1 & 1 & 0 & 1\ \\
\ 0 & 0 & 0 & 0 & 1 & 1 \
\end{matrix}\right), \quad
C=\left(\begin{matrix}
\ 1 & 3 &\ \ 5 & 1 & -2 & \ \ 2\ \\
\ 1 & 1 & -3 & 3 & \ \ 1 & -1 \\
\ 1 & 3 & \ \ 1 & 3 & -1 & \ \ 1
\end{matrix}\right).$$

Here $P_1=P_2=P_3=Q$ which is a prism depicted in \rf{prism}. We label the vertices of $Q$ using the order of the columns in $\bar A$.
\begin{figure}[h]
\includegraphics[scale=.3]{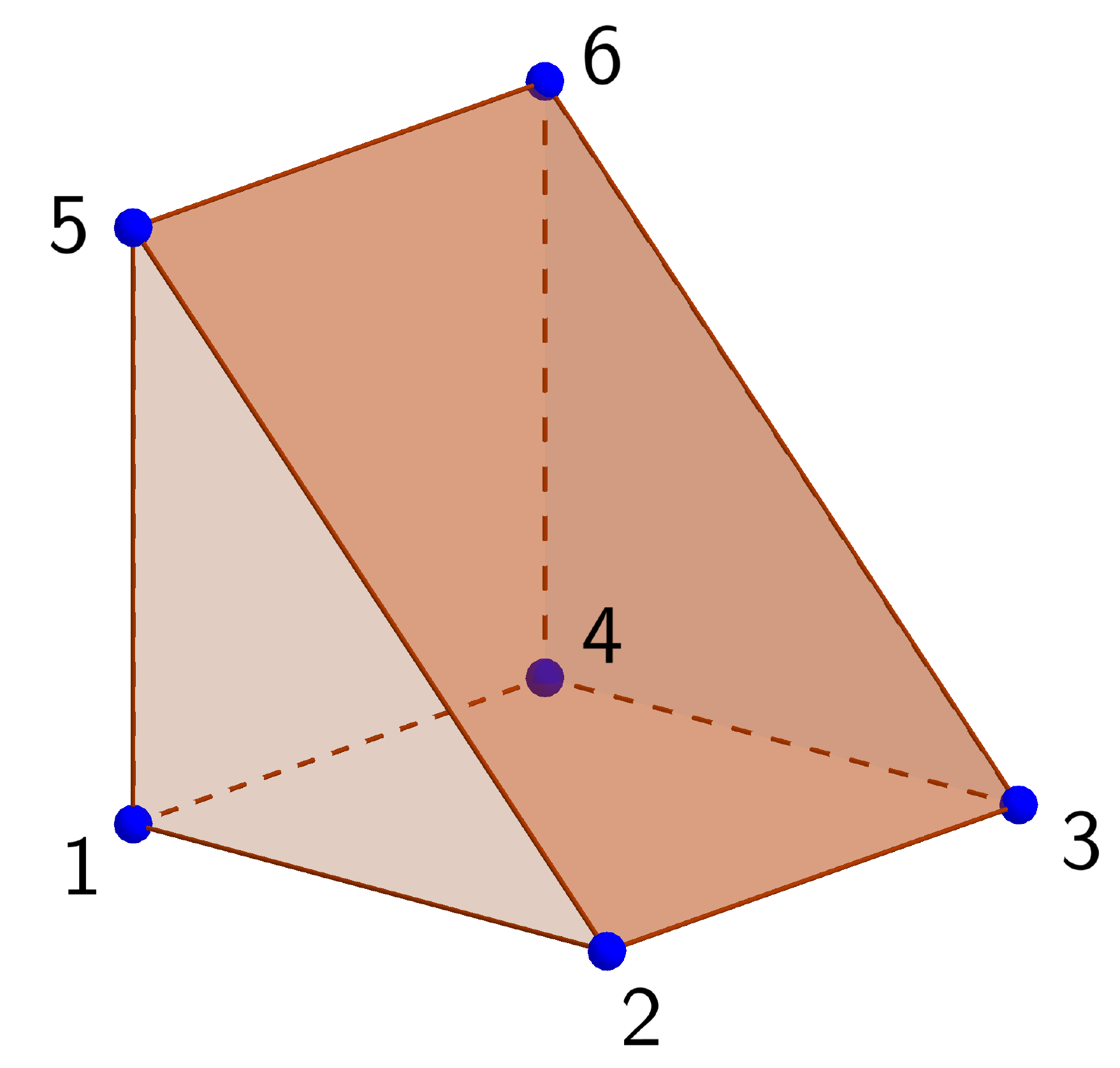}
\caption{The Newton polytope of the system in \rex{prism}.}
\label{F:prism}
\end{figure}
The submatrix of $C$ corresponding to the edge $F$ labeled $\{5,6\}$ has rank $1$ which is less than $\dim F +1$. Therefore the associated system has less than $3!\Vol_3(Q)=3$ isolated solutions in $(\C^*)^3$. (In fact, it has two solutions.) In particular, this is a degenerate system.
\end{Ex}
In the following very particular situation, the rank condition \eqref{E:rank} in Theorem \ref{T:Ber} implies the non-degeneracy of the system.

\begin{Rem} Assume that $P_1=P_2=\cdots=P_n=Q$ with $\dim Q=n$ and any proper face $F$ of $Q$ is a simplex which intersects $\cA$ only at its vertices.
Assume furthermore that $\text{\rm rank} \, C_{\cF} \geq \text{\rm rank} \, \bar{A}_{\cF}$ for any proper face $F$ of $Q$.
Then \eqref{E:system} is non-degenerate and thus has precisely $n!\,\Vol_n(Q)$ solutions in $(\K^*)^n$ counted with multiplicity according to Theorem \ref{T:BKK}.
Indeed, if $F$ is a proper face of $Q$, then the corresponding restricted system has total support $F \cap \cA$. If this restricted system is consistent, then
there is a non-zero vector in the kernel of $C_{\cF}$ and thus $\text{\rm rank} \, C_{\cF} < |F \cap \cA|=1+\dim F$ which gives a contradiction.
\end{Rem}

\providecommand{\bysame}{\leavevmode\hbox to3em{\hrulefill}\thinspace}
\providecommand{\MR}{\relax\ifhmode\unskip\space\fi MR }
\providecommand{\MRhref}[2]{%
  \href{http://www.ams.org/mathscinet-getitem?mr=#1}{#2}
}
\providecommand{\href}[2]{#2}

\end{document}